\documentclass[letterpaper, 10 pt, conference]{ieeeconf}  

\IEEEoverridecommandlockouts                              
\usepackage{cite}
\usepackage{url}
\usepackage{amsmath,amssymb,amsfonts}
\usepackage{multirow}
\usepackage{algorithmic}
\usepackage{graphicx}
\usepackage{textcomp}
\usepackage{epsfig}
\usepackage{algorithmic}
\usepackage{epstopdf}
\usepackage[ruled]{algorithm}
\usepackage{enumerate}
\usepackage{float}
\usepackage{graphics,graphicx}
\usepackage{subfigure}
\usepackage{booktabs}  

\usepackage{array}
\usepackage{mathtools}
\usepackage[justification=centering]{caption}
\allowdisplaybreaks[4]
\usepackage{url}
\DeclareGraphicsExtensions{.png}
\graphicspath{{./Pics/}}
\usepackage{color}
\usepackage{xspace}

\newcommand{\mb}{\mathbf}

\renewcommand{\Re}{{\mathbb{R}}}

\newtheorem{lemma}{Lemma}
\newtheorem{theorem}{Theorem}

\newtheorem{assumption}{Assumption}
\newtheorem{proposition}{Proposition}

\title{\LARGE \bf
Robust Distributed Nonconvex Optimization Enabling Communication Acceleration and Privacy Protection
}

\author{Zichong Ou and Jie Lu
\thanks{Z. Ou, J. Lu are with the School of Information Science and Technology, Shanghaitech University, 201210 Shanghai, China. J. Lu is also with the Shanghai Engineering Research Center of Energy Efficient and Custon AI IC, 201210 Shanghai, China. Email: {\tt\small ouzch, lujie@shanghaitech.edu.cn}.}
}

\begin{document}

\maketitle
\thispagestyle{empty}
\pagestyle{empty}

\begin{abstract} 
	This paper addresses a distributed nonconvex optimization problem over multi-agent networks, where each agent exchanges its local information solely with its neighbors. Given that most existing distributed nonconvex optimization algorithms are susceptible to information leakage during inter-agent communications, we propose a Robust Proximal Primal-dual algorithm, referred to as RPP, to enhance the security of information transmission. In contrast to many existing approaches that directly transmit local variables throughout the network, we introduce carefully designed random noises to obfuscate sensitive local information. This not only preserves privacy but also demonstrates the noise robustness of our proposed algorithm. We establish a sublinear rate at which RPP converges to a stationary solution. Moreover, by incorporating Chebyshev acceleration, an accelerated variant of RPP is developed and achieves the optimal communication complexity bound for the algorithms that allow for exchanging local decisions at each iteration. The superior convergence performance of RPP is validated through a few numerical experiments, which also indicate that, within an appropriate range, the introduced perturbations do not impede the convergence speed of RPP.
\end{abstract}

\section{INTRODUCTION}
Nonconvex optimization problems are ubiquitous in real-world scenarios, including optimal power flow in electric power systems \cite{molzahn2017survey}, network resource allocation \cite{tychogiorgos2011new} and wireless resource management \cite{lee2019deep}. The past decades have witnessed an explosive growth in data volume, which has spurred the development of distributed approaches. The field of distributed convex optimization is well-developed; however, nonconvex problems introduce significant new challenges for designing distributed algorithms.

In recent years, a few distributed nonconvex optimization algorithms have emerged. To handle nonconvexity, the Augmented Lagrangian (AL) function is widely employed in contemporary works \cite{NEXT2016,sun2016distributed,sun2018distributed,hong2017prox,alghunaim2022unified,sun2019distributed,mancino2023decentralized,yi2022sublinear,hong2016convergence}, which approximate nonconvex objectives via strongly convex surrogates and achieve sublinear convergence to stationary points. Among these, \cite{NEXT2016,sun2016distributed} use gradient tracking to boost convergence, while \cite{sun2018distributed,hong2017prox,alghunaim2022unified,sun2019distributed,mancino2023decentralized,yi2022sublinear,hong2016convergence} adopt primal-dual updates to strengthen consensus violation penalties—a core requirement for distributed optimization, which relies on neighbor communication. Most of the aforementioned algorithms \cite{ProxDGD2018,sun2018distributed,hong2017prox,alghunaim2022unified,yi2022sublinear,hong2016convergence} conduct one communication step and one computation step in each iteration, which becomes a bottleneck when the communication network is sparse due to the unbalanced convergence speed to reach consensus and stationarity. To tackle this limitation, \cite{xu2020accelerated,sun2019distributed,mancino2023decentralized} incorporate the well-known Chebyshev acceleration \cite{auzinger2011iterative} to optimize communication procedure. Additionally, \cite{sun2019distributed} proposes xFILTER—an optimal algorithm achieving the lower bound of communication complexity for distributed first-order methods—which is also attained by the ADAPD-OG-MC framework \cite{mancino2023decentralized}.

However, the aforementioned algorithms \cite{ProxDGD2018,NEXT2016,sun2016distributed,sun2018distributed,hong2017prox,alghunaim2022unified,sun2019distributed,mancino2023decentralized,yi2022sublinear,hong2016convergence} require nodes to directly transmit critical local information (e.g., decision variables or gradients) to neighbors, risking information leakage.
To address this privacy concern, several methods have been proposed. Some approaches \cite{gade2018privacy,wang2023decentralized,boyue2025TSP} employ additive perturbations to obfuscate shared gradients, \cite{lou2017privacy} estimates local variables by projecting the weighted average of local variables onto a certain set and \cite{wang2022decentralized} multiplies an additional stochastic mixing coefficient to local gradients when sharing information with neighbors. However, these strategies harm convergence: methods in \cite{wang2023decentralized,gade2018privacy,lou2017privacy} only guarantee asymptotic convergence (no explicit rates) due to diminishing stepsizes, while \cite{boyue2025TSP,wang2022decentralized} merely achieves convergence to a neighborhood of the stationary solutions.


In this paper, we propose a \emph{\underline{R}obust \underline{P}roximal \underline{P}rimal-dual algorithm (RPP)} for privacy-preserving distributed nonconvex optimization. Each node’s primal update minimizes a linearized Augmented Lagrangian function with an added proximal term, enabling multi-communication per iteration. We further incorporate a momentum-like mechanism into the dual update to boost convergence, and design an encryption strategy that injects well-designed random perturbations into shared messages to obfuscate local private variables. Our main contributions are as follows:
\begin{itemize}
	\item[1)] RPP, enhanced with our encryption strategy, is shown to converge to stationarity at a sublinear rate of $\mathcal{O}(\bar{M}\gamma^2/T)$, where $\bar{M}$ is the Lipschitz constant of the gradients of local functions and $\gamma\geq 1$ denotes the eigengap of the graph Laplacian matrix. As is shown in Table~\ref{tab:convergence rate}, RPP's convergence rate outperforms those of the methods without privacy protection (e.g., \cite{hong2017prox,yi2022sublinear}).
	\item[2)] We embed the Chebyshev acceleration in RPP, referred to as RPP-CA. It achieves the same state-of-the-art best convergence rate of $\mathcal{O}(\bar{M}/T)$ and the \emph{optimal communication complexity bound} of $\mathcal{O}(\bar{M}\sqrt{\gamma}/\epsilon)$ (for reaching an $\epsilon$-stationary error) as the methods in \cite{sun2019distributed,mancino2023decentralized}. This communication complexity bound is tailored to the class of algorithms where each node only communicates its local primal decision variables with its neighbors. 
    Notably, RPP-CA transmits variables of the same dimension per communication round as the methods in \cite{sun2019distributed,mancino2023decentralized}, while also ensuring data security without sacrificing efficiency.
	\item[3)] Our proposed encryption strategy preserves convergence rate and accuracy without increasing communication costs or extra assumptions. 
	In contrast, \cite{gade2018privacy} relies heavily on a client-server model—requiring servers to access all client information and causing additional communication.
    Additionally, \cite{lou2017privacy,wang2023decentralized} require bounded gradients of objective functions, while \cite{wang2022decentralized} demands Lipschitz continuity of local functions' Hessian matrices.
    \item[(4)] The numerical experiment against several benchmarks (e.g., \cite{hong2017prox,alghunaim2022unified,sun2019distributed,mancino2023decentralized,yi2022sublinear}) confirms the superior performance of our proposed algorithms with respect to iterations and communication rounds, which also demonstrates that our privacy mechanism does not impede the convergence speed when it is appropriately designed.
	\begin{table}[t]
    \centering
    \begin{tabular}{ccc}
        \toprule
        Method & Conv. rate & Comm. compl. \\
        \midrule
        L-ADMM \cite{yi2022sublinear} & $\mathcal{O}(\bar{M}^2\gamma^3/T)$ & $\mathcal{O}(\bar{M}^2\gamma^3/\epsilon)$ \\
        Prox-PDA \cite{hong2017prox} & $\mathcal{O}(\bar{M}^2\gamma^2/T)$ & $\mathcal{O}(\bar{M}^2\gamma^2/\epsilon)$\\
        ADAPD-OG \cite{mancino2023decentralized} & $\mathcal{O}(\bar{M}\gamma^2/T)$ & $\mathcal{O}(\bar{M}\gamma^2/\epsilon)$ \\
        ADAPD-OG-MC & $\mathcal{O}(\bar{M}/T)$ & $\mathcal{O}(\bar{M}\sqrt{\gamma}/\epsilon)$ \\
        xFILTER \cite{sun2019distributed} & $\mathcal{O}(\bar{M}/T)$ & $\mathcal{O}(\bar{M}\sqrt{\gamma}/\epsilon)$ \\
        \textbf{RPP} & $\mathcal{O}(\bar{M}\gamma^2/T)$ & $\mathcal{O}(\bar{M}{\gamma}^2/\epsilon)$\\
        \textbf{RPP-CA} & $\mathcal{O}(\bar{M}/T)$ & $\mathcal{O}(\bar{M}\sqrt{\gamma}/\epsilon)$\\
        \bottomrule
    \end{tabular}
    \caption{Convergence rates and communication complexity bounds of related works.}
    \label{tab:convergence rate}
\end{table}
	
\end{itemize}

The paper is structured as follows: Section~\ref{section problem form} introduces the nonconvex optimization problem, while Section~\ref{section algorithm design} elaborates on the development of RPP. Section~\ref{convergence analysis} states the convergence results for the proposed framework. Section~\ref{section simulation} compares the numerical results of RPP and several alternative methods, and Section~\ref{section conclusion} concludes the paper.

\emph{Notations:} For any differentiable function $f$, we denote its gradient by $\nabla f$. We represent the null space of a matrix argument as $\operatorname{Null}(\cdot)$. We denote the $n$-dimensional column all-one (all-zero) vector and identity (zero) matrix as $\mathbf{1}_n$ ($\mathbf{0}_n$) and $\mathbf{I}_n$ ($\mathbf{O}_n$), respectively, and the set of $n$-dimensional real symmetric matrices as $\mathbb{S}^n$.
In addition, we use $\langle \cdot,\cdot \rangle$, $\otimes$, and $\|\cdot\|$ for the Euclidean inner product, Kronecker product, and $\ell_2$ norm, respectively. For $\mathbf{A},\mathbf{B}\in \mathbb{R}^{d\times d}$, $\mathbf{A}\succ \mathbf{B}$ means $\mathbf{A}-\mathbf{B}$ is positive definite, and $\mathbf{A} \succeq \mathbf{B}$ means positive semi-definite.
For a matrix $\mathbf{A}$, we denote its $i$-th largest eigenvalue as $\lambda_i^\mathbf{A}$ and its Moore-Penrose inverse as $\mathbf{A}^{\dagger}$. For symmetric $\mathbf{A}\in \mathbb{R}^{d\times d}$ with $\mathbf{A} \succeq \mathbf{O}_d$ and $\mathbf{x}\in \mathbb{R}^d$, the weighted norm squared $\|\mathbf{x}\|^2_\mathbf{A}$ is defined as $\mathbf{x}^{\mathsf{T}}\mathbf{A}\mathbf{x}$.

\section{PROBLEM FORMULATION}\label{section problem form}

Let $\mathcal{G} = (\mathcal{V},\mathcal{E})$ be a connected, undirected graph, where $\mathcal{V} = \{1, \dots, N\}$ is the vertex set composed of $N$ nodes and the edge set $\mathcal{E} \subseteq \{\{i,j\}|i,j \in \mathcal{V},i \neq j\}$ describes the underlying interactions among the nodes. Through the network $\mathcal{G}$, each node $i$ exchanges information exclusively with neighbors in $\mathcal{N}_i=\{j\in\mathcal{V}:\{i,j\}\in\mathcal{E}\}$. All the nodes collaboratively solve the optimization problem:
\begin{equation}
\min _{x \in \mathbb{R}^{d}} f(x)= \sum_{i=1}^{N} f_{i}(x), \label{p1}
\end{equation}
where $f_i : \mathbb{R}^d \rightarrow \mathbb{R}$ is the local objective owned by node $i$. 

Next, We impose some assumptions on problem~\eqref{p1}:
\begin{assumption}\label{assumption smooth}
Each local objective $f_i\!\!:\!\!\Re^{d}\!\rightarrow\!\Re$ is differentiable and ${M}_i$-smooth, i.e., there exists ${M}_i> 0$ such that
\vspace{-0.3cm}
\begin{equation}
\|\nabla f_i(x) - \nabla f_i(y)\|\leq {M}_i\|x-y\|,\quad\forall x,y\in \Re^{d}.\label{smooth}
\end{equation}
\end{assumption}
\vspace{0.1cm}
\begin{assumption}\label{a1}
The function $f(x)$ is lower bounded over $x\in\mathbb{R}^d$, i.e., $
	f(x) \geq \inf_{x} f(x) > -\infty.$
\end{assumption}

Assumptions~\ref{assumption smooth} and~\ref{a1} are commonly adopted in existing works on distributed nonconvex optimization \cite{alghunaim2022unified, sun2018distributed, sun2019distributed,mancino2023decentralized, hong2016convergence,yi2022sublinear}. 

To address problem~\eqref{p1} over the graph $\mathcal{G}$, we assign to each node $i \in \mathcal{V}$ a local state variable $x_i \in \mathbb{R}^d$ as the estimate of the global decision $x \in \mathbb{R}^d$ in problem~\eqref{p1}, and define
\vspace{-0.1cm}
\begin{equation}
\tilde{f}(\mathbf{x}) := \sum _{i \in \mathcal{V}} f_i(x_i),\quad \mb{x} = (x_1^{\mathsf{T}}, \dots, x_N^{\mathsf{T}})^{\mathsf{T}} \in \Re ^{Nd}.\vspace{-0.1cm}
\end{equation}
As is demonstrated in \cite{Wu2020AUA}, problem~\eqref{p1} can be transformed into the following equivalent form:
\vspace{-0.1cm}
\begin{equation}
    \underset{\mathbf{x} \in \mathbb{R}^{N d}}{\operatorname{minimize}}\quad  \tilde{f}(\mb{x}) \qquad \text{s.t.}\quad \mb{L}^{\frac{1}{2}} \mb{x}=0, \label{p2}\vspace{-0.1cm}
\end{equation}
where $\mb{L}\in \mathbb{S}^{Nd}$ satisfies the following assumption.
\begin{assumption}\label{assumption L}
	The matrix $\mb{L}\in \mathbb{S}^{Nd}$ satisfies:
	\begin{itemize}
		\item[(\romannumeral1)] $\mb{L}$ is symmetric and positive semi-definite.
		\item[(\romannumeral2)] $\operatorname{Null}(\mb{L}) = \mathcal{S} := \{\mb{x} \in \Re ^{Nd}|x_1 = \cdots = x_N\}$.
	\end{itemize}
\end{assumption}

Assumption~\ref{assumption L} is commonly employed in distributed implementation \cite{alghunaim2022unified,mancino2023decentralized,yi2022sublinear}. Note that problem~\eqref{p1} and \eqref{p2} share the same optimal value. Clearly, under Assumption~\ref{assumption smooth}, $\tilde{f}$ is $\bar{M}-$smooth with $\bar{M}=\max_{i\in\mathcal
V}M_i$.

\section{ALGORITHM DEVELOPMENT}\label{section algorithm design}

This section develops a distributed algorithm for the optimization problem in \eqref{p2} to facilitate efficient network-wide information propagation while safeguarding node privacy.

To address the nonconvexity of the objectives, we employ the Augmented Lagrangian (AL) function: $\operatorname{AL}(\mb{x},\mb{v})=\tilde{f}(\mb{x})+(\mb{v})^{\mathsf{T}} \mb{L}^{\frac{1}{2}}\mb{x}+\frac{\rho}{2}\|\mb{x}\|^2_\mb{L}$, where $\mb{v}=(v_1^{\mathsf{T}}, \dots, v_N^{\mathsf{T}})^{\mathsf{T}}\in \Re ^{Nd}$ denotes the Lagrangian multiplier and $\rho>0$ is the penalty parameter. To address problem~\eqref{p2}, we employ the AL function and implement the following primal-dual framework: Starting from any $\mb{x}^0,\mb{v}^0,\hat{\mb{v}}^0\in \mathbb{R}^{Nd}$, for any $k \geq 0$,
\begin{align}
	\label{argmin approximation primal} \notag \mb{x}^{k+1}=&\underset{\mathbf{x} \in \mathbb{R}^{N d}}{\arg \min}\;\tilde{f}(\mb{x}^k) + \langle \nabla \tilde{f}(\mb{x}^k), \mb{x}-\mb{x}^k \rangle +\langle\mb{v}^{k}, \mb{L}^{\frac{1}{2}} \mb{x}\rangle\\ 
    &\qquad\qquad+\frac{\rho}{2}\|\mb{x}\|_{\mb{L}}^{2}+ \frac{1}{2} \|\mb{x}-\mb{x}^k\|^2_{\mb{B}}, \\ 
	\label{argmin approximation dual 1} \hat{\mb{v}}^{k+1}=&\hat{\mb{v}}^{k}+\rho \mb{L}^{\frac{1}{2}} \mb{x}^{k+1},\\
	\label{argmin approximation dual 2} \mb{v}^{k+1}=&\hat{\mb{v}}^{k+1}+\eta(\hat{\mb{v}}^{k+1}-\hat{\mb{v}}^{k}),
\end{align}
where $\eta \in \mathbb{R}$, and $\hat{\mb{v}}=(\hat{v}_1^{\mathsf{T}}, \dots, \hat{v}_N^{\mathsf{T}})^{\mathsf{T}}\in \Re ^{Nd}$ is a dual variable. The updates \eqref{argmin approximation primal}--\eqref{argmin approximation dual 2} can be interpreted as follows.
\begin{itemize}
    \item \textbf{Primal Update} \eqref{argmin approximation primal}: Linearize $\tilde{f}(\mb{x})$ at $\mb{x}^k$ as $\tilde{f}(\mb{x}^k)+\langle\nabla \tilde{f}(\mb{x}^k),\mb{x}-\mb{x}^k\rangle$ and incorporate a proximal term $\frac{1}{2}\|\mb{x}-\mb{x}^{k}\|^2_{\mb{B}}$ with $\mb{B}\in \mathbb{S}^{Nd}$ into the AL function.
    \item \textbf{Dual Ascent} \eqref{argmin approximation dual 1}: Execute a dual ascent step with the corresponding ``dual gradient" evaluated by the constraint residual at $\mb{x}^{k+1}$ generated by minimizing the AL-like function in \eqref{argmin approximation primal}.
    \item \textbf{Dual Acceleration} \eqref{argmin approximation dual 2}: Apply a momentum-like step to modify the dual ascent direction and boost convergence.
\end{itemize}

It is important to note that the primal update \eqref{argmin approximation primal} requires $\mb{B} + \rho \mb{L} \succ \mb{O}_{Nd}$ to ensure that $\mb{x}^{k+1}$ in \eqref{argmin approximation primal} is well-posed and admits a unique solution. Then, applying the first-order optimality condition, \eqref{argmin approximation primal} can be expressed as:
\vspace{-0.1cm}
\begin{equation}
\label{first-order opt} \nabla \tilde{f}(\mb{x}^k) + \mb{L}^{\frac{1}{2}} \mb{v}^k + \rho \mb{L} \mb{x}^{k+1}+ \mb{B}(\mb{x}^{k+1} - \mb{x}^k)=0.\vspace{-0.1cm}
\end{equation}

With the definition of $\mb{G}:= (\mb{B} + \rho \mb{L})^{-1}$ and \eqref{first-order opt}, we rewrite \eqref{argmin approximation primal} as
\vspace{-0.1cm}
\begin{equation}
	\mb{x}^{k+1} = \mb{x}^k - \mb{G}(\nabla \tilde{f}(\mb{x}^k) + \mb{L}^{\frac{1}{2}} \mb{v}^k + \rho \mb{L} \mb{x}^k). \label{x^k+1 vk}\vspace{-0.1cm}
\end{equation}

Note that \eqref{x^k+1 vk} and \eqref{argmin approximation dual 1} involve weight matrices $\mathbf{L}^{\frac{1}{2}}$ and $\mathbf{L}$—not distributively executable in a network. To allow information propagation across the network, we next describe how the algorithm is implemented in a distributed manner.

\subsection{Distributed Implementation}
We develop our distributed optimization algorithm based on \eqref{x^k+1 vk}, \eqref{argmin approximation dual 1} and \eqref{argmin approximation dual 2}. By introducing the variables $\mb{d}^k=\big((d_1^k)^{\mathsf{T}},\ldots,(d_N^k)^{\mathsf{T}}\big)^{\mathsf{T}}$ and $\hat{\mb{d}}^k=\big((\hat{{d}}_1^k)^{\mathsf{T}},\ldots,(\hat{{d}}_N^k)^{\mathsf{T}}\big)^{\mathsf{T}}$, we apply the following variable transformation:
\vspace{-0.1cm}
\begin{equation}
\mb{d}^k = \frac{1}{\rho}(\mb{L}^{\frac{1}{2}})^{\dagger} \mb{v}^k, \quad \hat{\mb{d}}^k = \frac{1}{\rho}(\mb{L}^{\frac{1}{2}})^{\dagger} \hat{\mb{v}}^k. \label{variable change}\vspace{-0.1cm}
\end{equation}
This requires $\mathbf{d}^k,\hat{\mathbf{d}}^k \in \mathcal{S}^{\perp}\, \forall k \geq 0$, where $\mathcal{S}^{\perp} := \{\mb{x}\in \Re ^{Nd}| \, x_1 + \cdots + x_N = \mb{0}\}$ is the orthogonal complement of $\mathcal{S}$, which can be simply guaranteed by
\vspace{-0.1cm}
\begin{equation}
\mb{d}^0, \hat{\mb{d}}^0 \in \mathcal{S}^{\perp} .\label{q0}\vspace{-0.1cm}
\end{equation}

By the variable substitution in \eqref{variable change}, the terms $\mb{L}^{\frac{1}{2}} \mb{v}^k + \rho \mb{L} \mb{x}^k$ in \eqref{x^k+1 vk} is substituted by $ \rho \mb{L}(\mb{d}^k+\mb{x}^k)$. Hence, from arbitrary $\mb{x}^0\in \mathbb{R}^{Nd}$ and $\mb{d}^0, \hat{\mb{d}}^0$ as in \eqref{q0}, for any $k \geq 0$,
\begin{align} 
\label{xk+1 original} &\mb{x}^{k+1} = \mb{x}^k - \mb{G} (\nabla \tilde{f}(\mb{x}^k)+\rho \mb{L}(\mb{x}^k+\mb{d}^k)),\\ 
\label{dtildek+1 original} &\hat{\mb{d}}^{k+1}=\hat{\mb{d}}^{k}+\mb{x}^{k+1},\\
\label{dk+1 original} &\mb{d}^{k+1}=\hat{\mb{d}}^{k+1}+\eta(\hat{\mb{d}}^{k+1}-\hat{\mb{d}}^{k}).\vspace{-0.1cm}
\end{align}

The iterations \eqref{q0}--\eqref{dk+1 original} are equivalent to \eqref{argmin approximation primal}--\eqref{argmin approximation dual 2}. Moreover, the property $\mb{B} + \rho \mb{L} \succ \mb{O}_{Nd}$ implies that the matrix $\mb{G}$ satisfies $\mb{G} \succ \mb{O}_{Nd}$. Note that $\mb{L}$ is a weight matrix responsible for local information exchange. To enable distributed implementation and enhance information propagation across the network, we design the matrix $\mb{G}$ as
\vspace{-0.1cm}
\begin{equation}
	\mb{G}=\alpha \mb{I}_{Nd}-\beta \mb{L}, \label{G=alpha I-beta L}\vspace{-0.1cm}
\end{equation}
where $\alpha>0, 0<\beta<\alpha/\lambda_1^{\mb{L}}$. The distributed implementation can be conducted by the weight matrix $\mb{L}$, defined as:
\vspace{-0.1cm}
\begin{equation}
\mb{L}= \mb{P} \otimes \mb{I}_d, \label{L=PI}\vspace{-0.1cm}
\end{equation}
where the symmetric matrix $\mb{P} \succeq \mb{O}_N$ has neighbor-sparse structures, i.e., its $(i,j)$-entries $p_{ij}$ vanish if $i\neq j$ and $\{i,j\}\notin\mathcal{E}$. Such a neighbor-sparse structure can be jointly determined by the nodes in a fully distributed manner \cite{Wu2020AUA}.

Subsequently, in the view of each node $i$, introducing $\mb{y}^k=\big((y_1^k)^{\mathsf{T}},\ldots,(y_N^k)^{\mathsf{T}}\big)^{\mathsf{T}}$, $\mb{z}^k=\big((z_1^k)^{\mathsf{T}},\ldots,(z_N^k)^{\mathsf{T}}\big)^{\mathsf{T}}$ and letting $\mb{y}^k=\mb{x}^k+\mb{d}^k$ and $\mb{z}^k = \nabla \tilde{f}(\mb{x}^k)+\rho \mb{L}\mb{y}^k$, the distributed implementation of \eqref{xk+1 original} can be written as
\vspace{-0.1cm}
\begin{align}
	y_i^k &=x_i^k+d_i^k,\label{yik}\\
	z_i^k &= \nabla f_i(x_i^k)+\rho \sum_{j \in \mathcal{N}_i \cup \{i\}} p_{ij} y_j^k,\label{zik}\\
	x_i^{k+1} &= x_i^k - \alpha z_i^k+\beta \sum_{j \in \mathcal{N}_i \cup \{i\}} p_{ij} z_j^k \label{xik+1}.\vspace{-0.1cm}
\end{align}
Updates \eqref{zik} and \eqref{xik+1} show that node $i$ transmits only $y_i^k$ and $z_i^k$ in iteration 
$k$, both containing local private information. 


\subsection{Privacy Protection Mechanism}
To further enhance data privacy, we incorporate perturbation variables $\mb{e}^{k}=\big((e_1^k)^{\mathsf{T}},\ldots,(e_N^k)^{\mathsf{T}}\big)^{\mathsf{T}}\in \mathbb{R}^{Nd}$ and $\mb{r}^{k}=\big((r_1^k)^{\mathsf{T}},\ldots,(r_N^k)^{\mathsf{T}}\big)^{\mathsf{T}}\in \mathbb{R}^{Nd}$ into $\mb{y}^k$ and $\mb{z}^k$, respectively, and rewrite \eqref{yik}--\eqref{zik} as 
\vspace{-0.1cm}
\begin{align}
    y_i^k &=x_i^k+d_i^k+e_i^k,\label{yik_new}\\
    z_i^k &= \nabla f_i(x_i^k)+\rho \sum_{j \in \mathcal{N}_i \cup \{i\}} p_{ij} y_j^k+r_i^k.\label{zik_new}\vspace{-0.1cm}
\end{align}
The variables $\mb{y}^k$ and $\mb{z}^k$ obscure local variables, preserving local privacy during the transmission of $\mb{x}^k+\mb{d}^k$ and $\nabla \tilde{f}(\mb{x}^k)+\rho \mb{L}\mb{y}^k$ in \eqref{xk+1 original}. Then, we rewrite \eqref{xk+1 original} as
\vspace{-0.1cm}
\begin{equation}\label{primal update final}
    \mb{x}^{k+1} = \mb{x}^k - \mb{G} (\nabla \tilde{f}(\mb{x}^k)+\mb{r}^k+\rho \mb{L}(\mb{x}^k+\mb{d}^k+\mb{e}^k)).\vspace{-0.1cm}
\end{equation}
Here, we impose the following assumption on $\mb{e}^{k}$ and $\mb{r}^k$.
\begin{assumption}\label{assumption ek+1}
    For $k\geq 1$ and some $\sigma_e,\sigma_r>0$, the perturbations variables ${e}_i^{k}$ and ${r}_i^{k}$ on nodes $i\in\mathcal{V}$ satisfy
    \vspace{-0.4cm}
    \begin{align}
	\label{ek+1 bound} \|{e}_i^{k}\|^2&\leq\sigma_e^2\|{x}_i^{k}-{x}_i^{k-1}\|^2,\\
	\label{ek+1 - ek bound}\|{e}_i^{k}-{e}_i^{k-1}\|^2&\leq \sigma_e^2\|{x}_i^{k}-{x}_i^{k-1}\|^2,\\
        \label{rk+1 bound} \|{r}_i^{k}\|^2&\leq\sigma_r^2\|{x}^{k}_i-{x}_i^{k-1}\|^2,\\
        \label{rk+1 - rk bound} \|{r}_i^{k}-{r}_i^{k-1}\|^2&\leq \sigma_r^2\|{x}_i^{k}-{x}_i^{k-1}\|^2.
    \end{align}
\end{assumption}

Assumption~\ref{assumption ek+1} is relatively restrictive than the perturbations in \cite{gade2018privacy,wang2023decentralized,boyue2025TSP,9354563}, yet it guarantees exact convergence at a high rate—unattainable in these works. It also requires diminishing noises, which may compromise privacy protection. This inherent limitation exists in privacy-preserving methods (e.g., \cite{wang2023decentralized,9354563}) and can be practically addressed by terminating iterations upon achieving the desired accuracy.

In each iteration, nodes conduct two encrypted data exchanges with neighbors. Using noises $\mathbf{e}^{k}$ and $\mb{r}^{k}$, they safeguard local information privacy during communication. Furthermore, as is analyzed in Section~\ref{convergence analysis}, our proposed algorithm is robust to perturbed transmission. We thus term it the \emph{\underline{R}obust \underline{P}roximal \underline{P}rimal-dual algorithm}, referred to as RPP. Its distributed implementation is detailed in Algorithm~\ref{algrithm DPP2}.

\floatname{algorithm}{Algorithm}
\begin{algorithm}[h]
        \renewcommand{\thealgorithm}{1}
	\caption{RPP}    
	\label{algrithm DPP2}               
	\begin{algorithmic}[1]
	  \STATE \textbf{Parameters:} $\rho,\alpha,\beta,\sigma_e,\sigma_r>0$, $\eta\in\Re$, $\mb{P}\succeq\mathbf{O}_N$.
		\STATE \textbf{Initialization:} Each node $i \in \mathcal{V}$ sets $x_i^{-1}=d_i^0=\hat{d}_i^0=x_i^0 =\mb{0}_d$.
		\FOR{$k \geq 0$} 
			\STATE Each node $i \in \mathcal{V}$ generates $e_i^k$ by \eqref{ek+1 bound}, \eqref{ek+1 - ek bound}, computes $y_i^k$ by \eqref{yik_new} and sends it to every neighbor $j \in \mathcal{N}_i$.
			\STATE Each node $i \in \mathcal{V}$ generates $r_i^k$ by \eqref{rk+1 bound}, \eqref{rk+1 - rk bound}, computes $z_i^k$ by \eqref{zik_new} and sends it to every neighbor $j \in \mathcal{N}_i$.
			\STATE Each node $i \in \mathcal{V}$ computes $x_i^{k+1}$ by \eqref{xik+1}.
			\STATE Each node $i \in \mathcal{V}$ computes $\hat{d}_i^{k+1}=\hat{d}_i^{k}+x_i^{k+1}$.
			\STATE Each node $i \in \mathcal{V}$ computes $d_i^{k+1}\!=\!\hat{d}_i^{k+1}\!+\!\eta(\hat{d}_i^{k+1}\!-\!\hat{d}_i^{k})$.
		\ENDFOR
	\end{algorithmic}
\end{algorithm}

\subsection{Design of communication acceleration scheme}\label{section DPP2-CA}

In this subsection, we attempt to accelerate the communication procedure via the design of weight matrix $\mb{L}$. 

Let $\kappa_{\mb{L}} := \lambda_1^{\mb{L}} / \lambda_{N-1}^{\mb{L}} \geq 1$ denote the eigengap of the matrix $\mb{L}$. From a graph connectivity perspective, $\mb{L}$ corresponds to a weighted interaction graph, and $\kappa_{\mb{L}}$ reflects the density of the graph. Note that the eigengap of a fully connected graph equals $1$. Intuitively, a smaller $\kappa_{\mb{L}}$ indicates stronger graph connectivity, which in turn facilitates information fusion. This observation is supported by Theorem~\ref{theorem UMAP-Opt sublinear} in Section~\ref{convergence analysis}, and is consistent with recent studies such as \cite{xu2020accelerated,sun2019distributed,mancino2023decentralized}.

Motivated by this interrelationship, by introducing a Laplacian matrix $\mb{H}=\mb{P}\otimes \mb{I}_d$ that satisfies Assumption~\ref{assumption L}, we redefine $\mb{L}$ from \eqref{L=PI} as a polynomial of $\mb{H}$, i.e.,
\begin{equation}
\mb{L} = P_{\tau}(\mb{H})/\lambda_1^{P_{\tau}(\mb{H})}, \label{L=PH}
\end{equation}
where $\tau\geq 1$ is the degree of polynomial and thus $\mb{L}$ satisfies Assumption~\ref{assumption L} with $\lambda_1^{\mb{L}}=1$. To mitigate the eigengap of $\mb{L}$, we incorporate the well-known Chebyshev acceleration technique \cite{auzinger2011iterative} into our distributed setting. This method allows us to generate an advantageous polynomial $P_{\tau}(\mb{H})$ for a fixed $\tau$. Practically, for any $\mb{s}\in \mathbb{R}^{Nd}$, the product $P_{\tau}(\mb{H})\mb{s}$ can be computed via local node interactions (See Oracle~\ref{oracle cheb}), enabling distributed computation of $\mb{x}^{k+1}$ according to \eqref{primal update final}.

\floatname{algorithm}{Oracle}
\begin{algorithm}[h] 
    \renewcommand{\thealgorithm}{$\mathcal{A}$}
    \caption{Chebyshev Acceleration} \label{oracle cheb}
    \begin{algorithmic}[1]
        \STATE \textbf{Input:} $\mb{s}=(s_1^{\mathsf{T}},\ldots,s_N^{\mathsf{T}})^{\mathsf{T}}\in \Re^{Nd}$, $\mb{P}\succeq\mathbf{O}_N$, $\tau=\lceil \sqrt{\kappa_{\mb{P}}} \rceil$, $c = \frac{\kappa_{\mb{P}}+1}{\kappa_{\mb{P}}-1}$ with $\kappa_{\mb{P}}=\lambda_1^{\mb{P}}/\lambda_{N-1}^{\mb{P}}$.\\
        \STATE \textbf{Procedure} $\text{CACC}(\mathbf{s},\mb{P},\tau)$ \\
        \STATE Each node $i \in \mathcal{V}$ computes $b^0 = 1$, $b^1 = c$.\\ 
	\STATE Each node $i \in \mathcal{V}$ maintains a variable $s_i^t$, sets $s_i^0 = s_i$ and sends it to every  neighbor $j \in \mathcal{N}_i$.\\
	\STATE Each node $i \in \mathcal{V}$ computes $s_i^1 = cs_i^0 - c\sum_{j\in {\mathcal{N}}_i \cup \{i\}}$ $ p_{ij} s_j^0$ and sends it to every  neighbor $j \in \mathcal{N}_i$.\\
        \FOR{$t=1:\tau -1$} 
        \STATE Each node $i \in \mathcal{V}$ computes $b^{t+1}= 2c b^t-b^{t-1}$.\\
        \STATE Each node $i \in \mathcal{V}$ computes $s_i^{t+1}= 2c s_i^t-s_i^{t-1}-2c\sum_{j\in {\mathcal{N}}_i \cup \{i\}} p_{ij} s_j^t$ and sends it to every neighbor $j \in \mathcal{N}_i$.\\
        \ENDFOR
        \STATE \textbf{Output:} Each node $i \in \mathcal{V}$ returns $s_i^0 - s_i^{\tau}/b^{\tau}$, so that $P_{\tau}(\mb{H})\mb{s}=\big((s_1^0 - s_1^{\tau}/b^{\tau})^{\mathsf{T}},\ldots,(s_N^0 - s_N^{\tau}/b^{\tau})^{\mathsf{T}}\big)^{\mathsf{T}}$.\\
        \STATE \textbf{End procedure}\\
    \end{algorithmic}
\end{algorithm}

By employing Oracle~\ref{oracle cheb} for communication, we update $\mb{z}^k$ in Line~5 of Algorithm~\ref{algrithm DPP2} for $\mb{z}^k=\nabla\tilde{f}(\mb{x}^k)+\rho\text{CACC}(\mb{y}^k,\mb{P},\tau)+\mb{r}^k$ and substitute \eqref{xik+1} in Line~6 for $\mb{x}^{k+1}=\mb{x}^k-\alpha \mb{z}^k+\beta \text{CACC}(\mb{z}^k,\mb{P},\tau)$. This leads to the development of the RPP method integrated with Chebyshev acceleration, referred to as RPP-CA. In each iteration, RPP-CA conducts $2\tau$ times communication rounds and has the potential to achieve high communication efficiency with an appropriate choice of $\tau$, which we will analyze in Section~\ref{convergence analysis}.

\section{CONVERGENCE ANALYSIS}\label{convergence analysis}
In this section, we provide the convergence analysis of RPP. For convenience, we let $\lambda_1^{\mb{L}}=1$, which can be easily extended to $\lambda_1^{\mb{L}}>0$. The following lemma connects the primal and dual variables, in which we define 
\begin{align}
	\label{wk+1} \mb{w}^{k+1} &:= (\mb{x}^{k+1}-\mb{x}^{k})-(\mb{x}^{k}-\mb{x}^{k-1}).
\end{align}

\begin{lemma} \label{lemma vk+1 - vk}
	Suppose Assumptions~\ref{assumption smooth}--\ref{assumption L} hold, and consider the use of \eqref{primal update final}, \eqref{dtildek+1 original} and \eqref{dk+1 original}. We also suppose $\mb{L}= \mb{P} \otimes \mb{I}_d$, $\mb{G}=\alpha \mb{I}_{Nd}-\beta \mb{L}$. Then, for all $k\geq 0$,
	\begin{align}
		&\rho\|\mb{x}^{k+1}\|^2_{\mb{L}}={\rho}\|\mb{d}^{k+1}-\mb{d}^k\|^2_{\mb{L}}\notag\\
		\leq& 5\kappa \big(d_1\|\mb{x}^k-\mb{x}^{k-1}\|^2_{(\mb{B})^{-1}}+\|\mb{w}^{k+1}\|^2_{\mb{B}}\big), \label{vk+1 - vk}
	\end{align}
	where $\kappa := \frac{\lambda_1^{\mb{B}}}{\rho \lambda_{N-1}^{\mb{L}}}$, $d_1 = \bar{M}^2+\rho^2(\eta^2+\sigma_e^2)+\sigma_r^2$.
\end{lemma}

By the variable change in \eqref{variable change}, we rewrite the AL function
\begin{equation}
	\operatorname{AL}(\mb{x},\mb{d})=\tilde{f}(\mb{x})+\langle \mb{x},\rho\mb{L}\mb{d} \rangle+\frac{\rho}{2}\|\mb{x}\|^2_\mb{L}. \label{AL d}
\end{equation}
Here, for convenience, we denote $\operatorname{AL}^k=\operatorname{AL}(\mb{x}^k,\mb{d}^k)$.
In the next lemma, we illustrates the dynamics of the AL function \eqref{AL d} for our proposed algorithm.

\begin{lemma} \label{lemma ALk+1 - ALk}
	Suppose all the conditions in Lemma~\ref{lemma vk+1 - vk} hold. For all $k\geq 0$,
	\begin{align}
		&\operatorname{AL}^{k+1}-\operatorname{AL}^{k} \leq  -\|\mb{x}^{k+1}-\mb{x}^{k}\|^2_{\mb{B}+\frac{1-\eta-\sigma_e}{2}\rho \mb{L}-\frac{\bar{M}+\sigma_r}{2}\mb{I}_{Nd}} \notag\\
		&+ \frac{5}{2}\kappa(2+\eta)\big(d_1\|\mb{x}^k-\mb{x}^{k-1}\|^2_{({\mb{B}})^{-1}}+\|\mb{w}^{k+1}\|^2_{\mb{B}}\big)\notag\\
		&-\frac{\eta\rho}{2}\|\mb{x}^k\|^2_{\mb{L}}+\frac{1}{2}\|\mb{x}^k-\mb{x}^{k-1}\|^2_{\sigma_r\mb{I}+\sigma_e\rho\mb{L}}. \label{ALk+1 - ALk lemma}
	\end{align}
\end{lemma}


Subsequently, we generate the following decreasing sequence (for some $c>0$):
\begin{align}
	&{P}^{k+1}\!=\!\operatorname{AL}^{k+1}\!+\|\mb{x}^{k+1}\!-\!\mb{x}^k\|^2_{(\frac{5}{2}(2+\eta)d_1\kappa)(\mb{B})^{-1}+\frac{1}{2}(\sigma_r\mb{I}+\sigma_e\rho\mb{L})}\notag\\
	&+\frac{{c}}{2}\big(\rho\|\mb{x}^{k+1}\|^2_\mb{L}+\|\mb{x}^{k+1}-\mb{x}^{k}\|^2_{\mb{B}+(\bar{M}+\sigma_r)\mb{I}_{Nd}+(|\eta|+\sigma_e)\rho \mb{L}}\big)\notag\\
	&-\frac{\eta\rho}{2}\|\mb{x}^{k+1}\|^2_{\mb{L}}. \label{tilde Pk+1}
\end{align}
In the next lemma, we show that, with proper parameters, the potential functions will decrease along iterations. 

\begin{lemma}\label{lemma tildePk+1 nonincreasing}
	Suppose all the conditions in Lemma~\ref{lemma vk+1 - vk} hold. For all $k\geq 0$, let the parameters of RPP satisfy:
	\begin{align}
		\label{eta, sigma} & |\eta|+2\sigma_e<\frac{1}{2}, \sigma_r\geq 0\\
		\label{define tilde c}& \kappa \leq cd_2/\lambda_{N-1}^{\mb{L}}, \text{ with } d_2=\frac{\lambda_{N-1}^{\mb{L}}}{6(2+[\eta]_+)},\\
		\label{B range}& \frac{1}{2}\mb{B}-\frac{1+2c}{2}(\bar{M}+2\sigma_r)\mb{I}_{Nd}-\frac{5d_1c}{12}\mb{B}^{-1}\succeq \mb{O}_{Nd}.
	\end{align}
	Then for all $k\geq 0$, we have
	\begin{equation}
		{P}^{k+1}\!-\!{P}^{k}\leq -\|\mb{x}^{k+1}\!\!-\!\mb{x}^k\|^2_{\frac{1}{4}(\mb{B}+\rho \mb{L})}-\frac{c}{12} \|\mb{w}^{k+1}\|^2_{\mb{B}}\leq 0. \label{tildePk+1 - tildePk lemma}
	\end{equation}
\end{lemma}


Next, we show the boundedness of $\{{P}^{k+1}\}$.
\begin{lemma}\label{lemma tildePk+1 geq f*}
	Suppose all conditions in Lemma~\ref{lemma vk+1 - vk} hold, and choose the parameters of \eqref{primal update final}, \eqref{dtildek+1 original} and \eqref{dk+1 original} by \eqref{eta, sigma}--\eqref{B range}. Let $\mb{x}^{0}=\hat{\mb{d}}^{0}=0$, for $k \geq 0$, the sequence $\{{P}^{k+1}\}$ satisfies
	\begin{equation}
		{P}^{k+1}\geq f^*, \,\,\, {P}^{1}\leq \tilde{f}(\mb{x}^1)+\frac{2+c}{(1+2c)(\bar{M}+2\sigma_r)}\|\nabla \tilde{f}(0)\|^2, \label{tildePk+1 geq f*}
	\end{equation}
	where $f^*$ is defined in Assumption~\ref{a1} and $\mb{x}^1 = -\mb{G}\nabla \tilde{f}(0)$.
\end{lemma}


Moreover, to satisfy \eqref{eta, sigma}--\eqref{B range}, we choose the algorithm parameters of \eqref{primal update final}, \eqref{dtildek+1 original} and \eqref{dk+1 original} in the following way:
\begin{align}
		\label{bounded eta, sigma} |\eta|<&\frac{1}{2} ,\quad \sigma_r>0,\quad 0\leq\sigma_e<\frac{1}{4}-\frac{1}{2}|\eta|,\\
		\label{bounded c}  c >& \frac{20\Delta^2}{3d_2^2}, \text{ with } \Delta >1,\\
		\label{bounded rho} \rho > & \frac{d_4+\sqrt{d_4^2+4d_3d_5}}{2d_3}, \text{ with } d_3=\frac{(d_2c+1)^2}{\Delta^2}-\frac{20c}{3},\notag\\
		 d_4 =& \frac{2(cd_2+1)(1+2c)(\bar{M}+2\sigma_r)}{\Delta}, \,\,\, d_5= \frac{10c\bar{M}^2}{3},\\
		\label{bounded alpha} \xi_1<&\frac{1}{\alpha}<\Delta \xi_1, \,\, \beta = \frac{\alpha}{2},\text{ with } \xi_1\!=\!\frac{1}{2}\Big( d_6\!+\!\sqrt{d_6^2\!+\!\frac{8d_1c}{3}} \Big),\notag\\
		  d_6&=(1+2c)(\bar{M}+2\sigma_r),
\end{align}
where $d_2$ is defined in \eqref{define tilde c}. We verify that the conditions above and the parameters in \eqref{eta, sigma}--\eqref{B range} are well-posed in the supplementary materials.

Subsequently, we state the convergence result of RPP in the following theorem.

\begin{theorem}\label{theorem UMAP-Opt sublinear}
	Suppose all the conditions and initialization in Lemma~\ref{lemma tildePk+1 geq f*} hold, and the parameters of iterations \eqref{primal update final}, \eqref{dtildek+1 original} and \eqref{dk+1 original} are chosen as \eqref{bounded eta, sigma}--\eqref{bounded alpha}, Then for $T>0$,	
	\begin{equation}
		\frac{1}{T}\sum_{k=1}^T(\frac{1}{N}\|\sum_{i=1}^N \nabla \tilde{f}_i(x_i^k)\|^2 + \rho\|\mb{x}^{k}\|^2_\mb{L})  \leq C_1 {C_2}/{T}, \label{epsilon leq C1C2/T}
	\end{equation}
	where $C_1:= \tilde{f}(\mb{x}^1)-f^*+\frac{2}{\bar{M}}\|\nabla \tilde{f}(0)\|^2$ and $C_2:= \frac{4}{\alpha}+8+\frac{1}{(2+[\eta]_+)c}(8+(5+\frac{18}{c})\frac{12}{c})$.
\end{theorem}


Theorem~\ref{theorem UMAP-Opt sublinear} shows that RPP converges to a stationary solution of problem~\eqref{p2} at a sublinear rate. Furthermore, when we set $\rho = \mathcal{O}(\kappa_\mathbf{L}\bar{M})$ according to \eqref{bounded rho}, we can specify $\frac{1}{T}\sum_{k=1}^T(\frac{1}{N}\|\sum_{i=1}^N \nabla \tilde{f}_i(x_i^k)\|^2 + \rho\|\mb{x}^{k}\|^2_\mb{L}) = \mathcal{O}(\kappa_{\mb{L}}^2\bar{M}/T)$, which is positively correlated with the eigengap $\kappa_{\mb{L}}$. This verifies our statement in Section~\ref{section DPP2-CA}. In the next lemma, we establish some results of Chebyshev acceleration.
\begin{lemma}\label{lemma Chebyshev}
	Suppose $\tau=\lceil \sqrt{\kappa_{\mb{P}}} \rceil$, $\mb{L}=P_{\tau}(\mb{H})/\lambda_1^{P_{\tau}(\mb{H})}$, we employ the Chebyshev polynomial to $P_{\tau}(\mb{H})$ as is shown in Oracle~\ref{oracle cheb}. Then, $\kappa_\mb{L}=\frac{{\lambda}_1^{P_{\tau}(\mb{H})}}{\lambda_{N-1}^{P_{\tau}(\mb{H})}}\leq (\frac{e^{\frac{1}{2}}+e^{-\frac{1}{2}}}{e^{\frac{1}{2}}-e^{-\frac{1}{2}}})^2.$
\end{lemma}

Lemma~\ref{lemma Chebyshev} implies that $\frac{1}{T}\sum_{k=1}^T(\frac{1}{N}\|\sum_{i=1}^N \nabla \tilde{f}_i(x_i^k)\|^2 + \rho\|\mb{x}^{k}\|^2_\mb{L}) = \mathcal{O}(\bar{M}/T)$ for RPP-CA. Considering that each iteration of RPP-CA conducts $2\lceil \sqrt{\kappa_{\mb{P}}} \rceil$ communication rounds, we next provide the communication complexity bound of RPP-CA.

\begin{proposition}\label{theorem optimal communication complexity bound}
	Suppose all the conditions in Theorem~\ref{theorem UMAP-Opt sublinear} hold and choose the parameters as \eqref{bounded eta, sigma}--\eqref{bounded alpha}. Additionally, for $k\geq 0$, let the iterations of \eqref{primal update final}, \eqref{dtildek+1 original} and \eqref{dk+1 original} be executed as is stated in Lemma~\ref{lemma Chebyshev}. Then, RPP-CA requires $T_c =\mathcal{O}(\frac{\bar{M}\sqrt{\kappa_{\mb{P}}}}{\epsilon})$ communication rounds to reach the $\epsilon$-stationary error.
\end{proposition}

The communication complexity bound above matches the optimal results in \cite{sun2019distributed,mancino2023decentralized}. However, these methods require each node $i$ to communicate its local variable $x_i$ to neighbors per iteration—risking information leakage. In contrast, our framework avoids exchanging $x_i$ entirely, propagating only encrypted information to significantly enhance security.

\section{NUMERICAL EXPERIMENT}\label{section simulation}
In this section, we evaluate the convergence performance of RPP and RPP-CA via a numerical example. 

We consider the distributed binary classification problem with nonconvex regularizers, adhering to Assumptions~\ref{assumption smooth}--\ref{a1}, formulated as \eqref{p1}, and characterized by
\begin{equation*}
	f_i(x)=\frac{1}{m}\sum_{s}^{m}\log(1+\exp(-y_{is}x^{\mathsf{T}}z_{is}))+\sum_{t=1}^{d}\frac{\lambda\mu([x]_t)^2}{1+\mu([x]_t)^2}.
\end{equation*}
Here, $m$ is the number of data samples of each node, $y_{is}\in \{-1,1\}$ and $z_{is}\in \mathbb{R}^d$ denote the label and the feature for the $s$-th data sample of node $i$, respectively. In the simulation, we set $N=50$, $d=10$, $m=200$, regularization parameters $\lambda=0.001$ and $\mu=1$, randomly generate $y_{is}$ and $z_{is}$ for each node $i$, and construct a geometric graph with $r=0.3$.

\begin{figure}[h]
    \centering
    \includegraphics[width=0.8\linewidth]{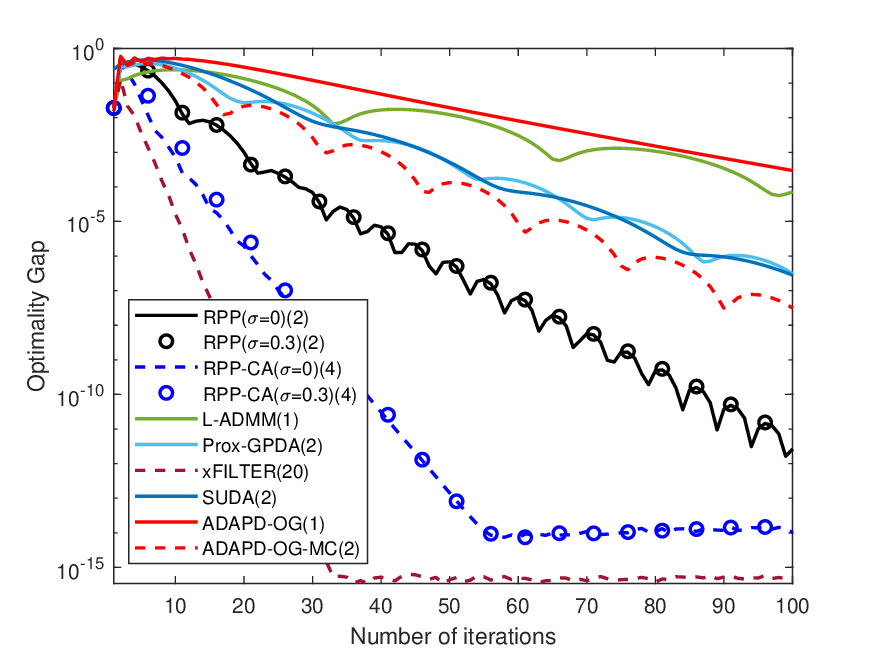}
    \caption{Convergence vs. iterations.}
    \label{iter}
\end{figure}

\begin{figure}[h]
    \centering
    \includegraphics[width=0.8\linewidth]{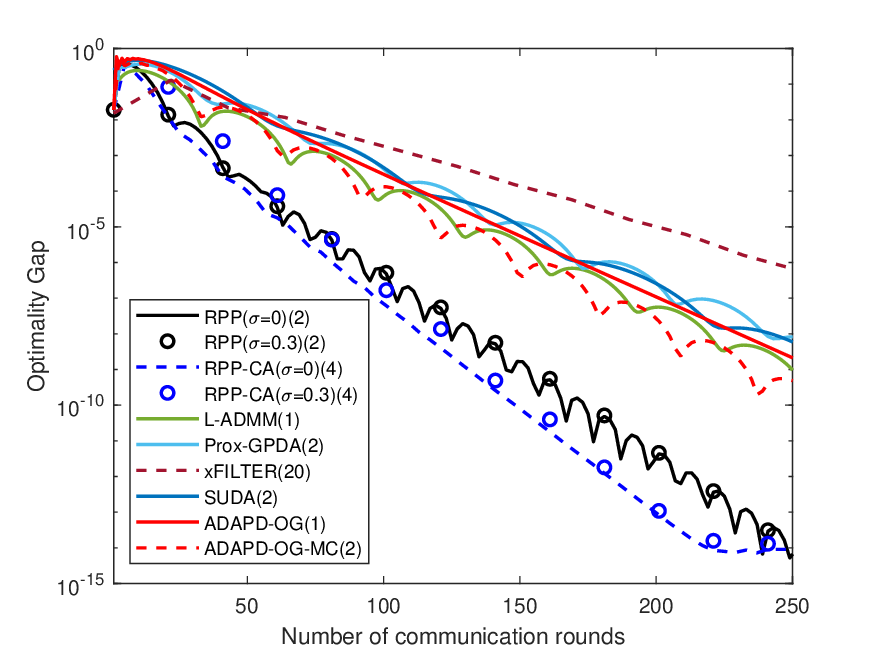}
    \caption{Convergence vs. communication rounds.}
    \label{comm}
\end{figure}

We compare RPP, RPP-CA with state-of-the-art distributed nonconvex optimization algorithms: L-ADMM \cite{yi2022sublinear}, SUDA \cite{alghunaim2022unified}, Prox-GPDA \cite{hong2017prox}, xFILTER \cite{sun2019distributed}, ADAPD-OG \cite{mancino2023decentralized} and its multi-communication variant ADAPD-OG-MC. Since RPP-CA, xFILTER and ADAPD-OG-MC conduct inner loops for Chebyshev acceleration per iteration, which increase communication cost, we compare their optimality gaps in terms of iterations (Fig.~\ref{iter}) and communication rounds (Fig.~\ref{comm}). 


For the implementation of these algorithms, we use the following parameter settings: In SUDA, we let $\mb{A}=\mb{I}-\mb{H},\mb{B}=\mb{H}^{1/2},\mb{C}=\mb{I}$; In RPP, we compare its performance for both $\sigma_r=\sigma_e=\sigma=0$ and $\sigma_r=\sigma_e=\sigma=0.3$, which satisfies Assumption~\ref{assumption ek+1}; In RPP-CA, we let $\tau=2$ for Chebyshev acceleration with the same two $\sigma$ configurations as RPP ($\sigma=0$ and $\sigma=0.3$); In ADAPD-OG-MC, we set its inner loop $R=2$. All algorithm parameters are manually optimized, and their performance is evaluated based on the optimality gap defined as $\|\nabla \tilde{f}(\mb{x})\|^2+\|\mb{H}^{\frac{1}{2}}\mb{x}\|^2.$

We summarize the numerical results as follows:
\begin{itemize}
    \item \textbf{Faster convergence speed and higher communication efficiency:} Fig.~\ref{iter} shows that our proposed RPP and RPP-CA converge faster than L-ADMM, Prox-GPDA, ADAPD-OG, and ADAPD-OG-MC, with the exception of xFILTER. The accelerated convergence of xFILTER, however, comes at the cost of communication overhead: xFILTER requires 20 internal communication rounds per iteration, leading to significantly lower communication efficiency, as is illustrated in Fig.~\ref{comm}. In contrast, RPP-CA requires the minimum number of communication rounds to reach stationarity, which reflects its highest communication efficiency among the compared algorithms.
    \item \textbf{Effectiveness of communication acceleration:} As is illustrated in Fig.~\ref{comm}, the comparative experiment between RPP-CA and RPP reveals that Chebyshev acceleration effectively enhances communication efficiency with an appropriately selected parameter $\tau$. This observation is further corroborated by the performance contrast between ADAPD-OG-MC (equipped with Chebyshev acceleration) and its parent algorithm ADAPD-OG. 
    \item \textbf{Robustness to perturbations:} The comparison of the perturbed (i.e., $\sigma=0.3$) and unperturbed (i.e., $\sigma=0$) versions of RPP reveals that both achieve stationarity in almost the same number of iterations and communication rounds. In addition, the accelerated variant, RPP-CA, yields analogous results. Collectively, these observations demonstrate that RPP and RPP-CA exhibit robustness to local perturbations within a reasonable range, while concurrently enhancing data privacy.
\end{itemize}

\section{CONCLUSIONS}\label{section conclusion}
We have proposed a robust proximal primal-dual algorithm for privacy protection, named as RPP, to minimize a distributed nonconvex, smooth optimization problem. We introduce an encryption strategy that incorporates additional perturbations to safeguard local information. RPP is shown to reach the stationarity at a sublinear rate and the encryption strategy does not obstruct the convergence. Furthermore, with Chebyshev acceleration, we demonstrate that the accelerated variant of RPP achieves an earlier derived communication complexity bound, which is optimal for the algorithm class where only local variables are transmitted through the network. The numerical result validates the high efficiency of our algorithms in terms of both iterations and communication rounds, while also demonstrating the robustness of RPP and RPP-CA against local communication perturbations.

\section{APPENDIX} \label{appendix}
\subsection{Proof of Lemma~\ref{lemma vk+1 - vk}}\label{proof lemma vk+1 - vk}
From \eqref{variable change}, \eqref{primal update final}, \eqref{dtildek+1 original}, \eqref{dk+1 original} and the first-order optimality condition \eqref{first-order opt}, for $k \geq 0$ we have
\begin{equation}
	\nabla \tilde{f}(\mb{x}^{k})+\mb{r}^k+\rho\mb{L}(\hat{\mb{d}}^{k+1}+\eta \mb{x}^{k}+\mb{e}^k)+\mb{B}(\mb{x}^{k+1}-\mb{x}^k)=0.  \label{optimal condition2}
\end{equation}
This implies
\begin{align}
	&\rho\mb{L}(\hat{\mb{d}}^{k+1}-\hat{\mb{d}}^k)=-(\nabla \tilde{f}(\mb{x}^k)-\nabla \tilde{f}(\mb{x}^{k-1}))-\mb{B}\mb{w}^{k+1}\notag\\
	&\!-\!\rho \mb{L}(\mb{e}^{k}\!-\!\mb{e}^{k-1})\!-\!\eta\rho \mb{L}(\mb{x}^k\!-\!\mb{x}^{k-1})\!-\!(\mb{r}^k\!-\!\mb{r}^{k-1}).\label{L1/2vk+1 - vk Bk=B}
\end{align}
Due to Assumption~\ref{assumption smooth}, \eqref{ek+1 - ek bound}, \eqref{rk+1 - rk bound} and $\lambda_1^{\mb{L}}=1$, we have
\begin{align*}
	&\rho\|\mb{x}^{k+1}\|^2_{\mb{L}}\!\overset{\eqref{dtildek+1 original}}{=}\!{\rho}\|\hat{\mb{d}}^{k+1}\!-\!\hat{\mb{d}}^k\|^2_\mb{L}
	\!\leq\!  \frac{1}{\rho \lambda_{N-1}^{\mb{L}}}\|{\rho}\mb{L}(\hat{\mb{d}}^{k+1}\!-\!\hat{\mb{d}}^k)\|^2\notag\\
    \leq& \frac{1}{\rho\lambda_{N-1}^{\mb{L}}}\|\nabla \tilde{f}(\mb{x}^k)-\nabla \tilde{f}(\mb{x}^{k-1})+\mb{B}\mb{w}^{k+1} \notag\\
	&+\rho \mb{L}(\mb{e}^{k}-\mb{e}^{k-1})+\eta\rho \mb{L}(\mb{x}^k-\mb{x}^{k-1})+(\mb{r}^k-\mb{r}^{k-1})\|^2  \notag \\
	\leq &\frac{5}{\rho\lambda_{N-1}^{\mb{L}}}\big(\|\nabla \tilde{f}(\mb{x}^k)-\nabla \tilde{f}(\mb{x}^{k-1})\|^2+\|\mb{w}^{k+1}\|^2_{\mb{B}^2}\notag\\
	&\!\!+\|\rho \mb{L}(\mb{e}^{k}\!-\!\mb{e}^{k-1})\|^2\!+\!\|\eta\rho \mb{L}(\mb{x}^k\!-\!\mb{x}^{k-1})\|^2\!+\!\|\mb{r}^k\!-\!\mb{r}^{k-1}\|^2\big) \notag\\
	\leq &\frac{5}{\rho\lambda_{N-1}^{\mb{L}}}\big(\bar{M}^2\|\mb{x}^k-\mb{x}^{k-1}\|^2+\|\mb{w}^{k+1}\|^2_{\mb{B}^2}\notag\\
	&+(\sigma_e^2\rho^2+\sigma_r^2)\|\mb{x}^k-\mb{x}^{k-1}\|^2+\eta^2\rho^2\| \mb{x}^k-\mb{x}^{k-1}\|^2\big)\notag\\
	\leq & 5\kappa\big(d_1\|\mb{x}^k-\mb{x}^{k-1}\|^2_{\mb{B}^{-1}}+\|\mb{w}^{k+1}\|^2_{\mb{B}}\big),
\end{align*}
where $\kappa := \frac{\lambda_1^{\mb{B}}}{\rho \lambda_{N-1}^{\mb{L}}}$, $d_1 = \bar{M}^2+\rho^2(\eta^2+\sigma_e^2)+\sigma_r^2$.

\subsection{Proof of Lemma~\ref{lemma ALk+1 - ALk}}\label{proof lemma ALk+1 - ALk}
Due to Assumption~\ref{assumption smooth}, \eqref{xk+1 original}, \eqref{AL d}, \eqref{ek+1 bound}, \eqref{rk+1 bound} and $\mb{G}=(\mb{B}+\rho \mb{L})^{-1}$, we have
\begin{align}\label{AL term1}
	& \operatorname{AL}(\mb{x}^{k+1},\mb{d}^{k})-\operatorname{AL}(\mb{x}^{k},\mb{d}^{k})\notag\\
	=&\tilde{f}(\mb{x}^{k+1})-\tilde{f}(\mb{x}^{k})+\langle \mb{x}^{k+1}-\mb{x}^k,\rho\mb{L}\mb{d}^k \rangle \notag\\
	&+\frac{\rho}{2}\|\mb{x}^{k+1}\|^2_\mb{L}-\frac{\rho}{2}\|\mb{x}^{k}\|^2_\mb{L}  \notag\\
	{\leq}& \langle \nabla \tilde{f}(\mb{x}^k)+\rho\mb{L}\mb{d}^k,\mb{x}^{k+1}-\mb{x}^k \rangle + \frac{\bar{M}}{2}\|\mb{x}^{k+1}-\mb{x}^k\|^2 \notag\\
    & + \frac{\rho}{2}\|\mb{x}^{k+1}\|^2_\mb{L} - \frac{\rho}{2}\|\mb{x}^{k}\|^2_\mb{L} \notag\\
	=& \langle \nabla \tilde{f}(\mb{x}^k)+\rho \mb{L}(\mb{d}^k+\mb{x}^k),\mb{x}^{k+1}-\mb{x}^k \rangle \!+\!\frac{\bar{M}}{2}\|\mb{x}^{k+1}\!-\!\mb{x}^k\|^2\notag\\
    &+ \frac{\rho}{2}\|\mb{x}^{k+1}\|^2_\mb{L} + \frac{\rho}{2}\|\mb{x}^{k}\|^2_\mb{L} -\langle \rho \mb{L} \mb{x}^k,\mb{x}^{k+1} \rangle \notag\\
	=& \!-\!\|\mb{x}^{k+1}\!-\!\mb{x}^k\|^2_{\mb{G}^{-1}-\frac{\bar{M}}{2}\mb{I}_{Nd}-\frac{1}{2}\rho \mb{L}}\!-\!\langle\mb{r}^k\!+\!\rho\mb{L}\mb{e}^k,\mb{x}^{k+1}\!-\!\mb{x}^k\rangle\notag\\
	\leq& \!-\!\|\mb{x}^{k+1}\!-\!\mb{x}^k\|^2_{\mb{B}\!+\!\frac{1\!-\!\sigma_e}{2}\rho \mb{L}\!-\!\frac{\bar{M}\!-\!\sigma_r}{2}\mb{I}_{Nd}}\!+\!\frac{1}{2}\|\mb{x}^k\!-\!\mb{x}^{k-1}\|^2_{\sigma_r\mb{I}\!+\!\sigma_e\rho\mb{L}}. 
\end{align}
From \eqref{dtildek+1 original} and \eqref{dk+1 original}, we have
\begin{align}\label{AL term2}
    &\operatorname{AL}(\mb{x}^{k+1},\mb{d}^{k+1})\!-\!\operatorname{AL}(\mb{x}^{k+1}\!,\mb{d}^k) \!=\! \langle\rho \mb{L}(\mb{d}^{k+1}\!\!-\!\mb{d}^k),\mb{x}^{k+1}\rangle\notag\\
    &\leq  \langle\rho \mb{L}\big(\hat{\mb{d}}^{k+1}-\hat{\mb{d}}^k+\eta(\mb{x}^{k+1}-\mb{x}^k)\big),\mb{x}^{k+1}\rangle \notag\\
    &=  (1+\frac{\eta}{2})\rho \|\mb{x}^{k+1}\|^2_{\mb{L}}\!+\!\frac{\eta \rho}{2} \|\mb{x}^{k+1}\!-\!\mb{x}^{k}\|^2_{\mb{L}}\!-\!\frac{\eta \rho}{2} \|\mb{x}^{k}\|^2_{\mb{L}}.
\end{align}
Combining \eqref{AL term1}, \eqref{AL term2} and using \eqref{vk+1 - vk}, we obtain \eqref{ALk+1 - ALk lemma}.

\subsection{Proof of Lemma~\ref{lemma tildePk+1 nonincreasing}}\label{proof lemma tildePk+1 nonincreasing}
From \eqref{L1/2vk+1 - vk Bk=B}, \eqref{dtildek+1 original}, \eqref{ek+1 - ek bound}, \eqref{rk+1 - rk bound}, \eqref{wk+1} and Assumption~\ref{assumption smooth},
\begin{align}
	0\!=& -\langle \rho \mb{L}\big( \hat{\mb{d}}^{k+1}\!-\!\hat{\mb{d}}^k\!+\! \eta(\mb{x}^k\!-\!\mb{x}^{k-1})\!+\!\mb{e}^{k}\!-\!\mb{e}^{k-1}\!+\! \mb{B}\mb{w}^{k+1} \big)\notag\\
	&+\nabla \tilde{f}(\mb{x}^k)-\nabla \tilde{f}(\mb{x}^{k-1})+\mb{r}^{k}-\mb{r}^{k-1}, \mb{x}^{k+1}-\mb{x}^{k} \rangle\notag\\
    \leq&-\frac{\rho}{2}\|\mb{x}^{k+1}\|^2_{\mb{L}}+\frac{\rho}{2}\|\mb{x}^{k}\|^2_{\mb{L}}-\frac{\rho}{2}\|\mb{x}^{k+1}-\mb{x}^{k}\|^2_{\mb{L}}\notag\\
    &+\frac{\rho|\eta|}{2}(\|\mb{x}^k-\mb{x}^{k-1}\|^2_{\mb{L}} + \|\mb{x}^{k+1}-\mb{x}^{k}\|^2_{\mb{L}} )\notag\\
    &+\frac{1}{2}(\|\mb{x}^k-\mb{x}^{k-1}\|^2_{\sigma_r\mb{I}+\sigma_e\rho\mb{L}} + \|\mb{x}^{k+1}-\mb{x}^{k}\|^2_{\sigma_r\mb{I}+\sigma_e\rho\mb{L}} )\notag\\
    &-\frac{1}{2}(\|\mb{x}^{k+1}-\mb{x}^{k}\|^2_{\mb{B}} - \|\mb{x}^{k}-\mb{x}^{k-1}\|^2_{\mb{B}} + \|\mb{w}^{k+1}\|^2_{\mb{B}})\notag\\
    &+\frac{\bar{M}}{2}(\|\mb{x}^k-\mb{x}^{k-1}\|^2+\|\mb{x}^{k+1}-\mb{x}^{k}\|^2)\notag\\
    \leq & \!-\!\big(\frac{1}{2}\|\mb{x}^{k+1}\!\!-\!\mb{x}^{k}\|^2_{\mb{B}-(\bar{M}\!+\!\sigma_r)\mb{I}_{Nd}+(1\!-\!|\eta|\!-\!\sigma_e)\rho \mb{L}}\! +\! \frac{\rho}{2}\|\mb{x}^{k+1}\|^2_{\mb{L}}\big)\notag\\
    & \frac{1}{2}\|\mb{x}^{k}-\mb{x}^{k-1}\|^2_{\mb{B}+(\bar{M}+\sigma_r)\mb{I}_{Nd}+(|\eta|+\sigma_e)\rho \mb{L}}+\frac{\rho}{2}\|\mb{x}^{k}\|^2_{\mb{L}}\notag\\
    &-\frac{1}{2}\|\mb{w}^{k+1}\|^2_{\mb{B}}.\label{term combine}
\end{align}
It follows from \eqref{tilde Pk+1}, \eqref{ALk+1 - ALk lemma}, \eqref{vk+1 - vk} and \eqref{term combine} that
\begin{align}
	&P^{k+1}-P^k \notag\\
	\underset{\eqref{tilde Pk+1}}{\overset{\eqref{ALk+1 - ALk lemma}}{\leq}} & -\|\mb{x}^{k+1}-\mb{x}^{k}\|^2_{\mb{B}+\frac{1-\eta-2\sigma_e}{2}\rho \mb{L}-\frac{\bar{M}+2\sigma_r}{2}\mb{I}_{Nd}} \notag\\
	&+ \frac{5}{2}\kappa(2+\eta)\big(d_1\|\mb{x}^k-\mb{x}^{k-1}\|^2_{({\mb{B}})^{-1}}+\|\mb{w}^{k+1}\|^2_{\mb{B}}\big)\notag\\
	&-\frac{\eta\rho}{2}\|\mb{x}^k\|^2_{\mb{L}}+\big(\frac{5}{2}(2+\eta)d_1\kappa\big)(\|\mb{x}^{k+1}-\mb{x}^k\|^2_{(\mb{B})^{-1}}\notag\\
	&-\|\mb{x}^{k}-\mb{x}^{k-1}\|^2_{(\mb{B})^{-1}}) -\frac{\eta\rho}{2}(\|\mb{x}^{k+1}\|^2_{\mb{L}}-\|\mb{x}^{k}\|^2_{\mb{L}})\notag\\
	&+\frac{c}{2}\Big(\rho \|\mb{x}^{k+1}\|^2_{\mb{L}}+\|\mb{x}^{k+1}-\mb{x}^{k}\|^2_{\mb{B}+(\bar{M}+\sigma_r)\mb{I}_{Nd}+(|\eta|+\sigma_e)\rho\mb{L}}\notag\\
	&-(\rho \|\mb{x}^{k}\|^2_{\mb{L}}+\|\mb{x}^{k}-\mb{x}^{k-1}\|^2_{\mb{B}+(\bar{M}+\sigma_r)\mb{I}_{Nd}+(|\eta|+\sigma_e)\rho\mb{L}})\Big)\notag\\
	\underset{\eqref{vk+1 - vk}}{\overset{\eqref{term combine}}{\leq}}&-\|\mb{x}^{k+1}-\mb{x}^{k}\|^2_{\mb{B}+\frac{1}{2}\big(1-\eta-2\sigma_e+c(1-2|\eta|-2\sigma_e)\big)\rho \mb{L}}\notag\\
	&+\|\mb{x}^{k+1}-\mb{x}^{k}\|^2_{-\frac{1+2c}{2}(\bar{M}+2\sigma_r)\mb{I}_{Nd}+\frac{5}{2}d_1\kappa(2+[\eta]_+)({\mb{B}})^{-1}}\notag\\
	&-\big(\frac{c}{2}-\frac{5}{2}\kappa(2+[\eta]_+)\big)\|\mb{w}^{k+1}\|^2_{\mb{B}}\notag\\
	\underset{\eqref{define tilde c}}{\overset{\eqref{eta, sigma}}{\leq}} &-\|\mb{x}^{k+1}-\mb{x}^{k}\|^2_{\mb{B}+\frac{1}{4}\rho \mb{L}-\frac{1+2c}{2}(\bar{M}+2\sigma_r)\mb{I}_{Nd}-\frac{5d_1c}{12}\mb{B}^{-1}}\notag\\
	&-\big(\frac{c}{2}-\frac{5}{2}\kappa(2+[\eta]_+)\big)\|\mb{w}^{k+1}\|^2_{\mb{B}}.
\end{align}
This, together with \eqref{B range} gives \eqref{tildePk+1 - tildePk lemma}.

\subsection{Proof of Lemma~\ref{lemma tildePk+1 geq f*}}\label{proof lemma tildePk+1 geq f*}
From \eqref{AL d}, \eqref{dtildek+1 original}, \eqref{dk+1 original}, $\lambda_1^{\mb{L}}=1$,
\begin{align}
	&\operatorname{AL}^{k+1} - \tilde{f}(\mb{x}^{k+1})=\langle \rho \mb{L}\mb{d}^{k+1},\mb{x}^{k+1}\rangle +\frac{\rho}{2}\|\mb{x}^{k+1}\|^2_{\mb{L}}\notag\\
	=& \langle \rho \mb{L}(\hat{\mb{d}}^{k+1}+\eta(\hat{\mb{d}}^{k+1}-\hat{\mb{d}}^{k})),\hat{\mb{d}}^{k+1}-\hat{\mb{d}}^{k}\rangle +\frac{\rho}{2}\|\mb{x}^{k+1}\|^2_{\mb{L}}\notag\\
	=& \frac{\rho}{2}(\|\hat{\mb{d}}^{k+1}\|^2_{\mb{L}}\!-\!\|\hat{\mb{d}}^{k}\|^2_{\mb{L}}\!+\!\|\hat{\mb{d}}^{k+1}\!-\!\hat{\mb{d}}^{k}\|^2_{\mb{L}})\!+\!\frac{1\!+\!2\eta}{2}\rho\|\mb{x}^{k+1}\|^2_{\mb{L}}\notag\\
	\geq & \frac{\rho}{2}(\|\hat{\mb{d}}^{k+1}\|^2_{\mb{L}}-\|\hat{\mb{d}}^{k}\|^2_{\mb{L}}) +(1+\eta)\rho\|\mb{x}^{k+1}\|^2_{\mb{L}}.
\end{align}

Define $\widehat{\operatorname{AL}}^{k+1}:=\operatorname{AL}^{k+1}-f^*$, $\hat{f}(\mb{x}):=\tilde{f}(\mb{x})-f^*\geq 0$, $\hat{P}^{k+1}:= {P}^{k+1}-f^*$. Summing over $k=0,\dots,T$, we obtain $\sum_{k=0}^T \widehat{\operatorname{AL}}^{k+1} \geq \frac{\rho}{2}(\|\hat{\mb{d}}^{T+1}\|^2_{\mb{L}}\!-\!\|\hat{\mb{d}}^{0}\|^2_{\mb{L}})+ \sum_{k=0}^T (\hat{f}(\mb{x}^{k+1})+(1+\eta)\rho\|\mb{x}^{k+1}\|^2_{\mb{L}}).$
By the initialization $\mb{x}^{0}=0$, $\hat{\mb{d}}^0=0$ and $|\eta|>\frac{1}{2}$, the above sum is lower bounded by zero. Therefore, the sum of $\hat{P}^{k+1}$ is also lower bounded by zero, i.e., $\sum_{k=0}^T \hat{P}^{k+1}\geq 0, \quad \forall T>0.$
Note that Lemma~\ref{lemma tildePk+1 nonincreasing} shows that $\hat{P}^{k+1}$ is nonincreasing under the parameter selections \eqref{eta, sigma} and \eqref{B range}. Therefore, we can conclude that
\begin{equation}
	\hat{P}^{k+1}\geq 0, \quad {P}^{k+1}\geq f^*, \quad \forall k \geq 0. \label{hat Pk+1 geq 0 tildePk+1 geq f*}
\end{equation}

Next, since $\mb{x}^{0}=0$, $\hat{\mb{d}}^{0}=0$, $\mb{e}^0=\mb{r}^0=0$ we have
\begin{align}
	{P}^1 =& \operatorname{AL}^{1}+\|\mb{x}^{1}\|^2_{(\frac{5}{2}(2+\eta)d_1\kappa)(\mb{B})^{-1}-\frac{\eta\rho}{2}\mb{L}}\notag\\
	&+\frac{{c}}{2}\|\mb{x}^{1}\|^2_{\mb{B}+(\bar{M}+2\sigma_r)\mb{I}_{Nd}+(1+|\eta|+\sigma_e)\rho \mb{L}}, \label{tildeP0} \\
	\operatorname{AL}^1=&\tilde{f}(\mb{x}^1)+\langle \rho\mb{L}\mb{d}^1,\mb{x}^1\rangle+\frac{\rho}{2}\|\mb{x}^1\|^2_\mb{L} \notag\\
	{=}& \tilde{f}(\mb{x}^1)+\langle\rho\mb{L}(\hat{\mb{d}}^0+(1+\eta)\mb{x}^1),\mb{x}^1\rangle +\frac{\rho}{2}\|\mb{x}^1\|^2_\mb{L}\notag\\
	\overset{\eqref{ek+1 bound}}{\leq}& \tilde{f}(\mb{x}^1)+(\frac{3}{2}+\eta)\rho\|\mb{x}^1\|^2_\mb{L}, \label{AL0 initialization}\\
	\mb{x}^1 =& -\mb{G}\nabla \tilde{f}(0) . \label{x0 initialization}
\end{align}
By incorporating \eqref{AL0 initialization} into \eqref{tildeP0} and using \eqref{eta, sigma}--\eqref{B range}, we have
\begin{align}
	{P}^1\leq &\tilde{f}(\mb{x}^1)+\|\mb{x}^1\|^2_{\frac{5}{2}(2+\eta)d_1\kappa\mb{B}^{-1}+(\frac{3}{2}+\frac{\eta}{2})\rho \mb{L}}\notag\\
	&+\frac{c}{2}\|\mb{x}^1\|^2_{\big((1+|\eta|+\sigma_e)\rho \mb{L}+\mb{B}+(\bar{M}+2\sigma_r)\mb{I}_{Nd}\big)}\notag\\
	\leq &\tilde{f}(\mb{x}^1)+\|\mb{x}^1\|^2_{(2+\frac{3}{4}c)\mb{G}^{-1}}
    \overset{\eqref{x0 initialization}}{\leq}  \tilde{f}(\mb{x}^1)+(2+c)\|\nabla\tilde{f}(0)\|^2_{\mb{G}}\notag\\
	\overset{\eqref{B range}}{\leq}& \tilde{f}(\mb{x}^1)+\frac{2+c}{(1+2c)(\bar{M}+2\sigma_r)}\|\nabla \tilde{f}(0)\|^2, \label{tildeP0 first}
\end{align}
where the last inequality holds since \eqref{B range} implies that $\mb{G}<\frac{1}{(1+2c)(\bar{M}+2\sigma_r)}\mb{I}_{Nd}$. Hence we obtain \eqref{tildePk+1 geq f*}.

\subsection{Proof of Theorem~\ref{theorem UMAP-Opt sublinear}}\label{proof theorem UMAP-Opt sublinear}
\textbf{(i)}First we illustrate the feasibility of the parameters in \eqref{eta, sigma}--\eqref{B range}. Obviously, \eqref{bounded eta, sigma} implies \eqref{eta, sigma} in Lemma~\ref{lemma tildePk+1 nonincreasing}. To see \eqref{B range}, it is sufficient to show
\begin{equation}\label{range lambdaNB}
	\frac{1}{2}\lambda_N^{\mb{B}}-\frac{1+2c}{2}(\bar{M}+2\sigma_r)\mb{I}_{Nd}- \frac{5d_1c}{12\lambda_N^{\mb{B}}}\geq 0.
\end{equation}
Since \eqref{bounded c}, we have $\rho < \frac{\beta \xi_1}{\alpha}=\frac{\xi_1}{2}$ for $\rho>0$. Hence, from $\lambda_i^{\mb{B}}=\frac{1}{\alpha-\beta \lambda_i^{\mb{L}}}-\rho \lambda_i^{\mb{L}}$, we obtain $\lambda_N^{\mb{B}}=\frac{1}{\alpha}$ and $\lambda_1^{\mb{B}}=\frac{1}{\alpha-\beta}-\rho=\frac{2}{\alpha}-\rho$. It follows from \eqref{bounded alpha} that $\lambda_N^{\mb{B}} > \frac{1}{2}\big( (1+2c)(\bar{M}+2\sigma_r)+\sqrt{(1+2c)^2(\bar{M}+\sigma_r)^2+\frac{10d_1c}{3}} \big)$, and thus \eqref{range lambdaNB} holds. Hence, we obtain \eqref{B range}.

Note that \eqref{define tilde c} is equivalent to $cd_2\rho > \lambda_1^{\mb{B}}=\frac{2}{\alpha}-\rho$ due to $\kappa = \frac{\lambda_1^{\mb{B}}}{\rho \lambda_{N-1}^{\mb{L}}},$. From \eqref{define tilde c} and \eqref{bounded rho}, we have $\rho > \frac{\Delta}{d_2c+1}\Big( d_6+\sqrt{d_6^2+\frac{10c}{3}(\bar{M}^2\!+\!2\rho^2+\sigma_r^2)} \Big),$ 
where $d_6=(1+2c)(\bar{M}+2\sigma_r)$. This implies that $\rho>\frac{\Delta\xi_1}{d_2c+1}$, since \eqref{bounded eta, sigma} infers $(\eta^2+\sigma_e^2)<2$. Together with $1/\alpha<\Delta\xi_1$ in \eqref{bounded alpha}, we obtain $\rho>\frac{1}{\alpha(d_2c+1)}$. Therefore, \eqref{define tilde c} is satisfied.

\textbf{(ii)}
Next, we analyze the convergence result of RPP. 

Here, we define $\mb{J}=\frac{1}{N}\mb{1}_N \mb{1}_N ^{\mathsf{T}} \otimes \mb{I}_d$. Thus, $\mb{J}$ and $\mb{B}$ are commutative in matrix multiplication due to Assumption~\ref{assumption L}. Then, We multiply both sides of \eqref{optimal condition2} by matrix $\mb{J}$, since $\mb{J}\mb{L}=0$ we obtain $\mb{J}\nabla \tilde{f}(\mb{x}^{k}) + \mb{J}\mb{B}(\mb{x}^{k+1}-\mb{x}^{k}),\forall k \geq 0.$
Next, from $\mb{J}\nabla \tilde{f}(\mb{x})=\mb{1}_N \otimes\big(\frac{1}{N} \sum_{i=1}^{N}\nabla\tilde{f}_i(x_i)\big)$, $\mb{J}\mb{J}=\mb{J}$, $\|\mb{J}\|=1$ and \eqref{tildePk+1 - tildePk lemma}, we have
\begin{align}
	&\frac{1}{N}\|\sum_{i=1}^N \nabla \tilde{f}_i(x_i^k)\|^2 \leq (\mb{x}^{k+1}-\mb{x}^{k})^{\mathsf{T}}\mb{B}\mb{J}\mb{J}\mb{B}(\mb{x}^{k+1}-\mb{x}^{k}) \notag\\
	\leq & \|\mb{x}^{k+1}-\mb{x}^{k}\|^2_{\mb{B}} \times \|\mb{J}\mb{B}\| \overset{\eqref{tildePk+1 - tildePk lemma}}{\leq}  \frac{4}{\alpha}(P^k-P^{k+1}). \label{bounded gradient norm}
\end{align}
From \eqref{range lambdaNB}, we have $\mb{B}\succeq \lambda_N^{\mb{B}} \mb{I} \succeq \frac{5d_1c}{6\lambda_N^{\mb{B}}}\mb{I} \succeq \frac{5d_1c}{6}\mb{B}^{-1}$. It then follows from \eqref{ALk+1 - ALk lemma} that
\begin{align}
	&\rho\|\mb{x}^{k+1}\|^2_\mb{L}\leq 5\kappa(d_1\|\mb{x}^k-\mb{x}^{k-1}\|^2_{\mb{B}^{-1}}+\|\mb{w}^{k+1}\|^2_{\mb{B}}) \notag\\
	\leq & 5\kappa(2d_1\|\mb{x}^{k+1}-\mb{x}^k\|^2_{\mb{B}^{-1}}+3d_1\|\mb{w}^{k+1}\|^2_{\mb{B}^{-1}}+\|\mb{w}^{k+1}\|^2_{\mb{B}})\notag \\
	\leq & 5\kappa(\frac{12}{5c}\|\mb{x}^{k+1}-\mb{x}^k\|^2_{\mb{B}}+(1+\frac{18}{5c})\|\mb{w}^{k+1}\|^2_{\mb{B}}). \notag
\end{align}
Combining this with \eqref{define tilde c}, \eqref{tildePk+1 - tildePk lemma}, we have
\begin{align}
	&\rho\|\mb{x}^{k}\|^2_\mb{L} \leq 2\rho\|\mb{x}^{k+1}\|^2_\mb{L}+2\rho\|\mb{x}^{k+1}-\mb{x}^k\|^2_\mb{L}\notag\\
	\overset{\eqref{define tilde c}}{\leq} & \frac{1}{6(2+[\eta]_+)}\big(\frac{12}{c}\|\mb{x}^{k+1}-\mb{x}^k\|^2_{\mb{B}}+(5+\frac{18}{c})\|\mb{w}^{k+1}\|^2_{\mb{B}}\big)\notag\\
    &+2\rho\|\mb{x}^{k+1}-\mb{x}^k\|^2_\mb{L}\notag\\
	\overset{\eqref{tildePk+1 - tildePk lemma}}{\leq}& (8+\frac{8+(5+\frac{18}{c})\frac{12}{c}}{(2+[\eta]_+)c})(P^k-P^{k+1}).\label{rho xk L}
\end{align}
Combining \eqref{bounded gradient norm} and \eqref{rho xk L}, and using \eqref{hat Pk+1 geq 0 tildePk+1 geq f*} and \eqref{tildeP0 first}, we have
\begin{align*}
	&\sum_{k=1}^T(\frac{1}{N}\|\sum_{i=1}^N \nabla \tilde{f}_i(x_i^k)\|^2 + \rho\|\mb{x}^{k}\|^2_\mb{L}) \notag\\
	\overset{\eqref{bounded gradient norm}\eqref{rho xk L}}{\leq} & \sum_{k=1}^{T}(\frac{4}{\alpha}+8+\frac{8+(5+\frac{18}{c})\frac{12}{c}}{(2+[\eta]_+)c})(P^k-P^{k+1})\notag\\
	\leq &C_2(P^1-P^{T+1})\overset{\eqref{hat Pk+1 geq 0 tildePk+1 geq f*}}{\leq} C_2(P^1-f^*)\notag\\
	\overset{\eqref{tildeP0 first}}{\leq}& C_2(\tilde{f}(\mb{x}^1)-f^*+\frac{2+c}{(1+2c)(\bar{M}+2\sigma_r)}\|\nabla \tilde{f}(0)\|^2). 
\end{align*}
Hence, we obtain \eqref{epsilon leq C1C2/T}.

\subsection{Proof of Lemma~\ref{lemma Chebyshev}}\label{proof lemma Chebyshev}
First, we state the properties of the Chebyshev polynomial that is inspired by Theorem 7 in \cite{xu2020accelerated}. Consider the normalized Laplacian $\mb{H}$ with a spectrum in $[1-c_1^{-1},1+c_1^{-1}]$ with $c_1=\frac{\kappa_{\mb{P}}+1}{\kappa_{\mb{P}}-1}$, the Chebyshev polynomial $P_{K}(x)=1-\frac{T_{\tau}(c_1(1-x))}{T_{K}(c_1)}$ introduced in \cite{auzinger2011iterative} is the solution of the following problem: $\min_{p\in \mathbb{P}^{K},p(0)=0} \max_{x\in[1-c_1^{-1},1+c_1^{-1}]} |p(x)-1|$.
Accordingly, we have $\max_{x\in[1-c_1^{-1},1+c_1^{-1}]} |P_K(x)-1| \leq 2\frac{c_0^K}{1+c_0^{2K}}$,
where $c_0=\frac{\sqrt{\kappa_{\mb{P}}}-1}{\sqrt{\kappa_{\mb{P}}}+1}$. By setting $K=\tau=\lceil \sqrt{\kappa_{\mb{P}}} \rceil$, we obtain $c_0^K=(\frac{\kappa_{\mb{P}}-1}{\kappa_{\mb{P}}+1})^{\lceil {\sqrt{\kappa_{\mb{P}}}} \rceil} \leq (\frac{\kappa_{\mb{P}}-1}{\kappa_{\mb{P}}+1})^{{\sqrt{\kappa_{\mb{P}}}}} \leq e^{-1}.$
It follows that $2\frac{c_0^K}{1+c_0^{2K}}=\frac{2}{c_0^K+(c_0^{K})^{-1}}\leq \frac{2}{e+e^{-1}}<1.$
Subsequently, we obtain the following property: $1-\frac{2}{e+e^{-1}}\leq 1-2\frac{c_0^K}{1+c_0^{2K}}\leq \lambda_{N-1}^{P_{\tau}(\mb{H})}\leq \lambda_{1}^{P_{\tau}(\mb{H})} 
    \leq 1+2\frac{c_0^K}{1+c_0^{2K}} \leq 1+ \frac{2}{e+e^{-1}}.$
Finally, we have derived Lemma~\ref{lemma Chebyshev} by $\kappa_\mb{L}={{\lambda}_1^{P_{\tau}(\mb{H})}}/{\lambda_{N-1}^{P_{\tau}(\mb{H})}}$.

\bibliographystyle{IEEEtran}
\bibliography{IEEEabrv,reference}

@ARTICLE{boyue2025TSP,
  author={Li, Boyue and Chi, Yuejie},
  journal={IEEE Journal of Selected Topics in Signal Processing}, 
  title={Convergence and Privacy of Decentralized Nonconvex Optimization With Gradient Clipping and Communication Compression}, 
  year={2025},
  volume={19},
  number={1},
  pages={273-282},
  doi={10.1109/JSTSP.2025.3526081}}

@article{Wu2020AUA,
  author={Wu, Xuyang and Lu, Jie},
  journal={IEEE Transactions on Automatic Control}, 
  title={A Unifying Approximate Method of Multipliers for Distributed Composite Optimization}, 
  year={2023},
  volume={68},
  number={4},
  pages={2154-2169},
  doi={10.1109/TAC.2022.3173171}
}

@inproceedings{sun2016distributed,
  author={Sun, Ying and Scutari, Gesualdo and Palomar, Daniel},
  booktitle={2016 50th Asilomar Conference on Signals, Systems and Computers}, 
  title={Distributed nonconvex multiagent optimization over time-varying networks}, 
  year={2016},
  volume={},
  number={},
  pages={788-794},
  doi={10.1109/ACSSC.2016.7869154}
}

@InProceedings{hong2017prox,
  title = 	 {Prox-{PDA}: The Proximal Primal-Dual Algorithm for Fast Distributed Nonconvex Optimization and Learning Over Networks},
  author =       {Mingyi Hong and Davood Hajinezhad and Ming-Min Zhao},
  booktitle = 	 {Proceedings of the 34th International Conference on Machine Learning},
  pages = 	 {1529--1538},
  year = 	 {2017},
  volume = 	 {70},
}

@article{sun2018distributed,
  title={Distributed Non-Convex First-Order Optimization and Information Processing: Lower Complexity Bounds and Rate Optimal Algorithms},
  author={Sun, Haoran and Hong, Mingyi},
  journal={arXiv preprint arXiv:1804.02729},
  year={2018}
}

@article{sun2019distributed,
  title={Distributed non-convex first-order optimization and information processing: Lower complexity bounds and rate optimal algorithms},
  author={Sun, Haoran and Hong, Mingyi},
  journal={IEEE Transactions on Signal processing},
  volume={67},
  number={22},
  pages={5912--5928},
  year={2019},
  publisher={IEEE}
}

@article{yi2022sublinear,
  author={Yi, Xinlei and Zhang, Shengjun and Yang, Tao and Chai, Tianyou and Johansson, Karl Henrik},
  journal={IEEE Transactions on Control of Network Systems}, 
  title={Sublinear and Linear Convergence of Modified ADMM for Distributed Nonconvex Optimization}, 
  year={2023},
  volume={10},
  number={1},
  pages={75-86},
  doi={10.1109/TCNS.2022.3186653}
}

@inproceedings{xu2020accelerated,
  title = 	 {Accelerated Primal-Dual Algorithms for Distributed Smooth Convex Optimization over Networks},
  author =       {Xu, Jinming and Tian, Ye and Sun, Ying and Scutari, Gesualdo},
  booktitle = 	 {Proceedings of the Twenty Third International Conference on Artificial Intelligence and Statistics},
  pages = 	 {2381--2391},
  year = 	 {2020},
  volume = 	 {108}
}

@article{auzinger2011iterative,
  title={Iterative solution of large linear systems},
  author={Auzinger, Winfried and Melenk, J},
  journal={Lecture notes, TU Wien},
  year={2011}
}

@article{mancino2023decentralized,
  title={A decentralized primal-dual framework for non-convex smooth consensus optimization},
  author={Mancino-Ball, Gabriel and Xu, Yangyang and Chen, Jie},
  journal={IEEE Transactions on Signal Processing},
  volume={71},
  pages={525--538},
  year={2023},
  publisher={IEEE}
}

@article{molzahn2017survey,
  title={A survey of distributed optimization and control algorithms for electric power systems},
  author={Molzahn, Daniel K and D{\"o}rfler, Florian and Sandberg, Henrik and Low, Steven H and Chakrabarti, Sambuddha and Baldick, Ross and Lavaei, Javad},
  journal={IEEE Transactions on Smart Grid},
  volume={8},
  number={6},
  pages={2941--2962},
  year={2017},
  publisher={IEEE}
}

@article{lee2019deep,
  title={Deep learning for distributed optimization: Applications to wireless resource management},
  author={Lee, Hoon and Lee, Sang Hyun and Quek, Tony QS},
  journal={IEEE Journal on Selected Areas in Communications},
  volume={37},
  number={10},
  pages={2251--2266},
  year={2019},
  publisher={IEEE}
}

@inproceedings{tychogiorgos2011new,
  title={A new distributed optimization framework for hybrid ad-hoc networks},
  author={Tychogiorgos, George and Gkelias, Athanasios and Leung, Kin K},
  booktitle={2011 IEEE GLOBECOM Workshops (GC Wkshps)},
  pages={293--297},
  year={2011},
  organization={IEEE}
}

@article{alghunaim2022unified,
  title={A unified and refined convergence analysis for non-convex decentralized learning},
  author={Alghunaim, Sulaiman A and Yuan, Kun},
  journal={IEEE Transactions on Signal Processing},
  volume={70},
  pages={3264--3279},
  year={2022},
  publisher={IEEE}
}

@article{hong2016convergence,
  title={Convergence analysis of alternating direction method of multipliers for a family of nonconvex problems},
  author={Hong, Mingyi and Luo, Zhi-Quan and Razaviyayn, Meisam},
  journal={SIAM Journal on Optimization},
  volume={26},
  number={1},
  pages={337--364},
  year={2016},
  publisher={SIAM}
}

@article{wang2023decentralized,
  title={Decentralized nonconvex optimization with guaranteed privacy and accuracy},
  author={Wang, Yongqiang and Ba{\c{s}}ar, Tamer},
  journal={Automatica},
  volume={150},
  pages={110858},
  year={2023},
  publisher={Elsevier}
}

@article{wang2022decentralized,
  title={Decentralized stochastic optimization with inherent privacy protection},
  author={Wang, Yongqiang and Poor, H Vincent},
  journal={IEEE Transactions on Automatic Control},
  volume={68},
  number={4},
  pages={2293--2308},
  year={2022},
  publisher={IEEE}
}

@article{lou2017privacy,
  title={Privacy preservation in distributed subgradient optimization algorithms},
  author={Lou, Youcheng and Yu, Lean and Wang, Shouyang and Yi, Peng},
  journal={IEEE transactions on cybernetics},
  volume={48},
  number={7},
  pages={2154--2165},
  year={2017},
  publisher={IEEE}
}

@inproceedings{gade2018privacy,
  title={Privacy-preserving distributed learning via obfuscated stochastic gradients},
  author={Gade, Shripad and Vaidya, Nitin H},
  booktitle={2018 IEEE Conference on Decision and Control (CDC)},
  pages={184--191},
  year={2018},
  organization={IEEE}
}

@ARTICLE{9354563,
  author={Ding, Tie and Zhu, Shanying and He, Jianping and Chen, Cailian and Guan, Xinping},
  journal={IEEE Transactions on Automatic Control}, 
  title={Differentially Private Distributed Optimization via State and Direction Perturbation in Multiagent Systems}, 
  year={2022},
  volume={67},
  number={2},
  pages={722-737},
  doi={10.1109/TAC.2021.3059427}}

@ARTICLE{ProxDGD2018,
  author={Zeng, Jinshan and Yin, Wotao},
  journal={IEEE Transactions on Signal Processing}, 
  title={On Nonconvex Decentralized Gradient Descent}, 
  year={2018},
  volume={66},
  number={11},
  pages={2834-2848},
  doi={10.1109/TSP.2018.2818081}}

@ARTICLE{NEXT2016,
  author={Lorenzo, Paolo Di and Scutari, Gesualdo},
  journal={IEEE Transactions on Signal and Information Processing over Networks}, 
  title={NEXT: In-Network Nonconvex Optimization}, 
  year={2016},
  volume={2},
  number={2},
  pages={120-136},
  doi={10.1109/TSIPN.2016.2524588}}

\end{document}